\documentclass[review]{article}

\usepackage{lineno,hyperref}
\modulolinenumbers[5]

\usepackage[latin9]{inputenc}
\usepackage{mathrsfs}
\usepackage{amsthm}
\usepackage{amsmath}
\usepackage{amssymb}
\usepackage{esint}

  \newcommand{\imp}{\mbox{$\ \Rightarrow\ $}}
\newcommand{\Imp}{\mbox{$\Longrightarrow$}}
\newcommand{\Iff}{\mbox{$\Longleftrightarrow$}}

\def\proof{{\bf Proof.}\ }
\def\ie{\emph{i.e.}}

\def\pes{\emph{e.g.}}

\def\D{{\mathcal D}}

\def\C{{\mathcal C}}
\def\F{{\mathcal F}}

\def\Gk{{\mathfrak G}}

\def\P{{\mathcal P}}
\def\Pf{{\mathcal P_{fin}}}

\def\U{{\mathcal U}}

\def\EE{{\mathbb E}}
\def\N{{\mathbb N}}

\def\R{{\,\mathbb R}}
\def\Z{{\mathbb Z}}
\def\Zak{{\mathbb Z}^{\kg}}
\def\Nak{{\mathbb N}^{\kg}}

\def\Zb_\kgk{{\mathbb Z}^{[\kg]}}

\def\WW{{\mathbb{W}}}

\def\WW{{\mathbb{W}}}

\def\LL{{\mathbb{L}}}

\def\Nk{{\mathfrak N}}

\def\ag{{\alpha}}
\def\dg{{\delta}}
\def\bg{{\beta}}
\def\eg{{\varepsilon}}

\def\kg{{\kappa}}
\def\og{{\omega}}

\def\cg{{\gamma}}
\def\sg{{\sigma}}
\def\Sg{{\Sigma}}
\def\thg{{\theta}}

\def\ah{\aleph}
\def\aho{\aleph_{0}}

\def\max{\mbox{\rm max}\;}

\def\fb{{\mathbf{f}}}
\def\Zb_\kg{{\mathbf{Z}}_\kg}

\def\xb{{\mathbf{x}}}

\def\Zb{{\mathbf{Z}}}
\def\ab{{\mathbf{a}}}
\def\bb{{\mathbf{b}}}

\def\nb{{\mathbf{n}}}

\def\tk{{\mathfrak{t}}}

\def\ak{{\mathfrak{a}}}
\def\ck{{\mathfrak{c}}}
\def\sk{{\mathfrak{s}}}
\def\mk{{\mathfrak{m}}}
\def\xk{{\mathfrak{x}}}
\def\yk{{\mathfrak{y}}}
\def\zk{{\mathfrak{z}}}

\def\nk{{\mathfrak{n}}}

\def\bk{{\mathfrak{b}}}
\def\pk{{\mathfrak{p}}}

\def\nk{{\mathfrak{n}}}

\def\Nk{{\mathfrak{N}}}

\def\*{\times}
\def\st{{}^{\star}\!\, }
\def\0{\emptyset}
\def\7{\setminus}
\def\_{\overline}
\def\eq{\simeq}
\def\<{\prec}

\def\ll{\preceq}
\def\lor{\vee}

\def\o+{\bigoplus}

\def\ult#1#2{^{#1}_{\; #2}}

\def\incl{\subseteq}

\def\linc{\supseteq}

\def\pincl{\subset}
\def\qincl{\sqsubseteq}
\def\qinclp{\sqsubset}
\def\lincq{\sqsupseteq}

\def\fin{\lhd}

\def\la{\langle}
\def\ra{\rangle}

\def\qed{\hfill $\Box$}

\def\CP{$\mathsf{CP} $}

\def\UP{$\mathsf{UP} $}

\def\diff{$\mathsf{(diff)} $}

\def\AP{\textsf{(AP)}}

\def\EP{\textsf{(EP)}}

\def\Ec{\textsf{(E5)}}

\def\CP{\textsf{(CP)}}
\def\PP{\textsf{(PP)}}
\def\UP{\textsf{(UP)}}
\def\WHP{\textsf{(WHP)}}

\def\part{\textsf{(part)}}

\def\oo{\textsf{(SupP)}}

\def\oob{\textsf{(SubP)}}

\def\diser{$(*)$}

\def\FIP{\textsf{FIP}}

\def\st{such that}

\def\Zl{{\mathcal Z}}
\def\Hl{{\mathcal H}}

\def\RA{$(\mathsf{RA}) $}

\def\SRA{$(\mathsf{SRA}) $}
\def\LA{$(\mathsf{LA}) $}
\def\CA{$(\mathsf{CA}) $}

\def\TI{(\mathsf{TI}) }

\def\PA{$(\mathsf{PA} )$}

\def\DSL{(\mathsf{DSL)}}
\def\DS{(\textsf{DSS)}}

\def\wh{\widehat}

\newtheorem{Theorem}{Theorem}[section]

\newtheorem{thm}[Theorem]{Theorem}

\newtheorem{lem}[Theorem]{Lemma}

\newtheorem{cor}[Theorem]{Corollary}

\theoremstyle{definition}
\newtheorem{defn}[Theorem]{Definition}

\newtheorem{rem}[Theorem]{Remark}

\begin{document}


\title{ Euclidean integers, Euclidean ultrafilters, and  Euclidea numerosities}

\author{Mauro Di Nasso, Marco Forti\footnote{
{Dipartimento di Matematica - Universit\`a di Pisa}}}
\maketitle



\begin{abstract}\label{abs}
We introduce axiomatically the ring $\Zb_\kg$ of the Euclidean integers, 
 that can be viewed as the ``integral part" of the field $\EE$ of Euclidean numbers of \cite{BFeu}, where the transfinite sum of ordinal indexed $\kg$-sequences of integers  is well defined.
  In particular 
  any ordinal  might be identified with 
  the transfiite sum of its characterstc funcion, preserving the so called natural operations. 

The ordered ring $\Zb_\kg$ may be obtained as  an ultrapower of $\Z$ modulo suitable ultrafilters, thus constituting a \emph{ring
of nonstandard integers.}
Most relevant is the 
    \emph{algebraic} characterization of the ordering: a Euclidean integer is  \emph{positive} if and only if it is  \emph{the transfinite sum of natural numbers.}
    This property requires the use of special ultrafilters called Euclidean, here intoduced to ths end.

    The ring $\Zb_\kg$ allows to  assign  a ``Euclidean"  size (\emph{numerosity}) 
to
 ``ordinal Punktmengen", \ie\ sets of tuples of ordinals, as 
the transfinite sum of their characteristic functions: so every set becomes equinumerous to a set of ordinals,
 the Cantorain defiitions of \emph{order, addition and multiplication} are maintained, while  the Euclidean principle ``the whole is 
greater than the part" is fulfilled.

%


\end{abstract}

%

\section*{Introduction}\label{intr}

The paper \cite{BFeu} introduces the ordered field $\EE$ of  the \emph{Euclidean numbers}, constituted by all the \emph{transfinite sums} of real numbers, of length less than
a fixed (strong) limit cardinal $\kg$. 
Similarly, we introduce here axiomatically, for any cardinal $\kg$, the ordered domain $\Zb_\kg$ of the \emph{Euclidean integers}, characterized as the collection of all the transfinite sums of ordinal indexed $\kg$-sequences of integers.
 Only \emph{ordinal indexed} sums are considered, so as to avoid the antinomies and
paradoxes that might affect the summing up of infinitely many numbers.
In particular,
the choice of ordinal numbers as indices seems particularly appropriate,
given their natural \emph{wellordered }structure, combined with the \emph{lattice }structure 
inherited from their finite subsets (see Subsection \ref{record}).

We call $\Zb_\kg$ the ring of the \textit{Euclidean integers}
because it arises in a ``numerosity
theory" of ordinal-labelled sets, whose main aim is to
save all the Euclidean common notions, including the fifth
    \emph{``the whole is greater
    than the part", }
 but still maintaining the Cantorian definitions of \emph{ordering, addition and multiplication} of sets.
Having at disposal transfinite sums of integers, a natural way of assigning numbers (\emph{numerosities}) to sets is 
to take the transfinite sums of their respective characteristic functions, which are 
Euclidean integers after an appropriate labelling of sets by ordinals. These numerosities  being nonnegative {Euclidean integers}, 
their arising arithmetic,  in contrast to the awkward  Cantorian cardinal arithmetic, shows the best algebraic properties, since these numerosities constitute the non-negative part of the ring $\Zb_\kg$,  a semiring of \emph{hypernatural numbers} 
 of Nonstandard Analysis (see \cite{BDNF8}).
                                
\medskip
The paper is organized as follows.\label{index}

 In Section \ref{euint} we introduce axiomatically the ring $\Zb_\kg$ of the Euclidean integers,  an ordered non-Archimedean
ring with a supplementary structure given
by the operation of \emph{transfinite sum} $\ \sum_\ag a_\ag$,
where  $\ \la a_{\ag}\mid\ag<\kg\ra\ $ is a $\kg$-sequence integers,
subject to four axioms stated in Subsection \ref{axiom}.
 Every Euclidean integer is obtained as a transfinite sum of ordinary integers, and
 more generally, any transfinite sum of integers is well defined. Moreover the
ring of the Euclidean integers is characterized, in Subsection \ref{limult},  as  an \emph{ultrapower of 
ordinal-indexed finite partial sums}, modulo suitable ultrafilters, hence $\Zb_\kg$ is a ring of \emph{hyperintegers}.

In Subsection \ref{embord} any  
ordinal 
 $\bg<\kg$ is associated to the transfinite
sum $\sum_\ag\chi_\bg(\ag)$ of its characteristic function $\chi_\bg(\ag)=\begin{cases}
  1    & \text{if}\ \ag<\bg, \\
  0   & \text{otherwise}.
\end{cases}$,
consistently with the so called \emph{natural ordinal operations} $\oplus$ and $\otimes$, so
the ring $\Zb_\kg$  might be considered as a sort of \emph{natural}
extension both of the ring of the ordinary integers $\Z$  and of the \emph{semiring} 
of ordinals $(\kg;\oplus,\otimes)$.
%
%
 
 Subsection \ref{cons} is dedicated to the main result of the paper, namely the existence of Euclidean ultrafilters, or equivalently of rings $\Zb_\kg$ of Euclidean integers satisfying the strong axiom \SRA, which allows for the most wanted identification of the non-negative part $\Zb_\kg^{\ge 0}$ with the set of all transfinite sums of \emph{natural numbers.}
 
  In Section \ref{eusiz} we specialize
the general Euclidean principles for dealing with the size of sets, and
 we explicit the axioms of our Euclidean theory, grounded on 
these   principles, together with appropriate properties of multiplication and (total) ordering  of \emph{numerosities}. In particular we deal with the proper superset property \oo, and with the equivalent difference property  \diff, which follow from the strong representation axiom \SRA, 

In Subsection \ref{numeu} Euclidean integers are assigned as numerosities to all point sets of finitely dimensional linear spaces over an ``ordinal line" in such a way that all properties corresponding to the Euclidean common notions are satisfied.
Moreover one obtains the supplementary benefits that every point set is equinumerous to a set of ordinals,  and conversely  that every nonnegative Euclidean integer $\xk$ is the numerosity of a set $X$ of ordinals, namely the transfinite sum of the characteristic function $\chi_X$.

A few final remarks and open questions can be found in Section
\ref{froq}.

\smallskip
In general, we refer to \cite{Je} for the set-theoretical notions and
facts used in this paper, and to \cite{CK} for definitions and facts
concerning ultrapowers, ultrafilters, and nonstandard models.

%
%
%
%
%
%

\section{The Euclidean integers}\label{euint}
We present here an axiomatization of the ordered domain $\Zb_\kg$ of the Euclidean integers, where all and only the transfinite sums of $\kg$-sequences of ordinary integers represent any element of $\Zb_\kg$. This axiomatization is independent of, but essentially the same as the Euclidean numbers of \cite{BFeu}.
The main new result is the possibility of postulating the stronger axiom \SRA, which allow to make the semiring $\Zb_\kg^{\ge 0}$ of the nonnegative Euclidean integers  coincide exactly with the set of all transfinite sums of natural numbers.

\subsection{Recalls on ordinals}\label{record}

The so called \emph{``natural"} product  and  sum of ordinals will be denoted by $\ag\otimes\bg$ and $\ag\oplus\bg$, 
respectively, whereas $\ag\bg$, $\ag+\bg$, and $\ag^\bg$  denote  the ordinary \emph{ordinal} operations. 

Recall that, given ordinals $\ag,\bg$, there exist uniquely determined 
ordinals $\cg\le\ag$ and $\dg<2^{\bg},$ such that
\[
\ag=\left(2^{\bg} \cg\right)+ \dg.
\]
%

Hence each ordinal has a unique \emph{base-$2$ normal form}
$$
\ag=\sum_{n=1}^{N}2^{\ag_{n}},\ 
where \ {i}<{j}\Rightarrow \ag
{}_{{i}}>\ag{}_{{j}},$$
and one has $\ag=\bigoplus_{n=1}^{N}2^{\ag_{n}}$,\ 
independently of the ordering of the exponents. 
%

In particular  $2^\og=\og$, and  the power $2^\ag=\og^\ag$ whenever $\ag= \og\ag$. It follows that the fixed points of the function $\ag\mapsto 2^\ag$ are $\og$ and and the so called $\eg$-numbers $\eg$ such that $\og^\eg=\eg$.


Recalling  the \emph{antilexicographic wellordering} of finite sets of ordinals
$$L < L'\ \  if\  and\ only\ if\ \ \max (L \bigtriangleup L' ) \in L',$$
 the $\ag$th finite set of ordinals, is
$$L_\ag =\{ \ag_1, ... ,\ag_n\}\ \  where \ \ \ag=\bigoplus_1^n 2^{\ag_i}.$$ 
In particular 
$$L_0=\0,\ \ \ L_{2^\ag}=\{\ag\},\ \ \ and\ \ \ L_{2^\ag+\bg}=\{\ag\}\cup L_\bg\ \ for\ all\ \bg<2^\ag.$$
So  $\Pf(\ag)$
can be naturally indexed by $2^\ag$.\
The correspondence $\alpha\mapsto L_\ag$  induces two restrictions $\qinclp,\
\fin $ of the ordinal ordering,  called  respectively \emph{formal inclusion},\footnote{
~The name formal inclusion should  also recall that the respective 
base-$2$ normal forms are indeed 
contained one inside of the other one.}
\ corresponding to ordinary set-inclusion between finite sets $\ \ag\qinclp\bg\iff L_\ag\pincl L_\bg$, and
 \emph{formal membership}   corresponding to ordinary membership 
 $\ \ag\fin\bg\iff\ag\in L_\bg$. Then we have:
\begin{enumerate}
   \item  
    Given $\ \ag=\o+_{i\in I}2^{\ag_{i}},\ \bg=\o+_{h\in H}2^{\ag_{h}},\ $
one has $\ \ag\sqsubset
\bg\ 
\Longleftrightarrow\ \, I\pincl H.\ $ and\ 
$\ \ag\fin\bg\ \Iff\ 2^\ag\qincl\bg\iff \exists h\in H \ \ag=\ag_h$.
\smallskip

 \item For $A\incl\kg$ let $\wh A=\ \{\ag\mid L_\ag\incl A \}$:
    then the map $\ell:\ag\mapsto L_\ag$ is a lattice isomorphism of
\  $(\wh A,\qinclp)$ onto $([A]^{<\og},\pincl)$. 
 \  In particular, for $|X|=\kg$,
 
 $~~~([X]^{<\og},\pincl)\ \cong\ ([\kg]^{<\og},\pincl)\ \cong (\kg, \sqsubset)\ \cong
  (\eg,\qinclp)$ for $\kg\le\eg=2^\eg<\kg^+$.
 
\item The following  properties hold:
\begin{itemize}
  \item $\,0\sqsubseteq \ag$
for all $\ag$, and $\ \ag\qinclp\bg\imp\ag<\bg$ ;
  \item  $|\{\bg\mid \bg\qincl \ag\}|=2^{|L_\ag|}$ and $ |\{\bg\mid \bg\fin \ag\}|={|L_\ag|}$
are \emph{finite} for all 
    $\ag$;

 \item for  all $\ag,\bg,\cg$ one has\ $\ \ \ag,\bg\fin\cg\ \ \Iff\ \ 2^\ag\oplus 2^\bg\qincl \cg;$
\item  for all $\ag$ and all $\bg,\cg,\xi<2^\thg$  one has the following useful criteria:

\begin{description}
    \item[\textbf{(C)}] 
$~~~
   ~~\,2^{\thg} \ag+\bg\,\qinclp 
    \dg\ \Iff \        2^{\thg} \ag,\, \bg\qinclp \dg\,$ 
    \item[\textbf{(D)}] 
 $~~~\dg=2^{\thg\cdot 2}\ag+2^\thg\xi+\xi\ \Imp\ \Big(\,2^{\thg} \cg+\bg\,\qinclp 
    \dg\ \Iff \         \cg,\, \bg\qinclp \dg\,\Big)$ 
\end{description}
\end{itemize}
  \item The \emph{cones }\ $C(\thg)=\{\ag\mid\thg\qinclp\ag\}$, for $\thg<\kg$,  
 generate a filter $\C_\kg$ on $\kg$: call \emph{fine} a filter $\F$ on $\kg$ that contains $\C_\kg$.
\item 
For $\eta<\kg$, let $\ D(\eta)=\{\dg=2^{\eta\cdot 2}\ag+2^\eta\xi+\xi\mid \xi<2^\eta,\,\ag<\kg\}$:\   
the family of sets\ 
$ D(\eta,\thg)=D(\eta)\cap C(\thg)$\ has
 the \FIP,  and generates a\emph{ fine filter }$\D_\kg$ on $\kg$: call \emph{superfine} a filter $\F$ on $\kg$ that contains $\D_\kg$. 
\end{enumerate}

\subsection{Axiomatic introduction of the Euclidean integers} \label{axiom}
From the algebraic point of view, the {Euclidean integers} are a \emph{non-Archimedean discretely 
ordered superring $\Zb_\kg$ of $\Z$}, with a supplementary structure, the \emph{Euclidean structure},
 introduced axiomatically
\emph{via} the operation of \emph{transfinite sum} 
$\Sg(\ab)=\sum_{\ag}a_{\ag},$
where $\ab=\la a_\ag\mid\ag<\kg\ra\in\Z^\kg$  is 
 any \emph{$\kg$-sequence} of integers 
 
 We don't need the full field $\EE$ of the Euclidean numbers, so we do not  ground on their general theory of \cite{BFeu}. We present here instead a simpler independent axiomatization of the ring of the \emph{Euclidean integers} $\Zb_\kg$, which is best suited for assigning Euclidean sizes to sets, in particular if integrated with the stronger axiom \SRA, whose consistency was not known before.
 
 Let $\Zb_\kg$ denote a \emph{discretely ordered commutative domain}, endowed  with the supplementary ``Euclidean" structure given by the \emph{transfinite sum}
  $$\sum_{\ag}a_{\ag}=\Sg(\ab),\   \textrm{ where}\ \ \la a_{\ag}\mid \ag<\kg\ra=\ab\in\Zak,$$
  
 Remark that we intend that any transfinite sum comprehends \emph{all summands} $a_\ag,\, \ag<\kg$. 
When needed, we restrict the sum to a subset $K\incl\kg$ by putting
$$\sum_{\ag\in K} a_\ag=\sum_{\ag} b_\ag,\ \ \textrm{with} \ \ b_\ag=a_\ag\cdot\chi_K(\ag), \ \textrm{and}\ \
\chi_K(\ag)=\begin{cases}
    1  & \text{if } \ag\in K,  \\
    0  & \text{otherwise}.
\end{cases}
$$  

 We make the natural assumption that  a \emph{transfinite sum} coincides with the 
\emph{ordinary
 sum} of the ring $\Zb_\kg$ when the number of \emph{non-zero summands} is \emph{finite}, and, similarly to \cite{BFeu},  we postulate for the ring ${\Zb_\kg}$ the following axioms.\footnote{~
For sake of clarity, we tend to denote general Euclidean numbers by  \emph{fractures}
 $\ak,\bk,\ck,{\sk},\tk,\xk,\yk,\zk$, and integers by \emph{latin}
letters $a,b,c,m,n,p,q,u,v,w,x,y,z$. The ordinal indices are denoted by
 greek letters ${\alpha},{\beta},{\gamma},{%
\delta},\eta$, with $\kg,\nu,\mu$ reserved for cardinals; $\kg$-sequences are denoted by boldface letters $\ab, \bb,\xb$, \ldots} 

First, an axiom  representing each Euclidean integer as a transfinite sum of (ordinary) integers
\begin{description}
\item[\RA ](\emph{\textbf{Representation Axiom}})\\ 
For any $ \xk$ in $\Zb_\kg$\ there\ exists\ 
$\
\mathbf{x}=\la x_{\ag}\mid \ag<\kg\ra\in\,\Z^{\kg}$  such\ that\  $\xk= \sum_{\ag} x_\ag.$
\end{description}

Next, a linearity axiom:
\begin{description}
\item[\LA] ({\emph{\textbf{Linearity Axiom}}) }\ The transfinite sum is 
${\,\mathbb{Z}}$-linear, \emph{i.e.}  
$$u\!\sum_{\ag}x_{\ag}+v\!\sum_{\ag}y_{\ag}=\sum_{\ag}(ux_{\ag}+vy_{\ag})\ \ \ for\ all\ \ u,v,
x_\ag,y_\ag\in{\Z}.$$
Then an axiom for \emph{comparing} transfinite sums:
\item[\CA] ({\emph{\textbf{Comparison Axiom}})} ${~}$
For all $\mathbf{x,y}\in{\,\Zak}$ 
\begin{equation*}
 ~\ \exists\,\thg\!<\!\kg\  \forall \, \dg\lincq \thg\, \Big( \ \sum_{\ag\qincl \dg}x_{\ag}\,\leq\,\sum_{\ag\qincl\dg
}y_{\ag}\,\Big)
   \ \Imp\ ~ \sum_{\ag}x_{\ag}\leq\sum_{\ag}y_{\ag}
\end{equation*}
\noindent(Remark that any sum $\sum_{\ag\qincl\dg}x_\ag$ is an \emph{ordinary finite
sum of integers}.)

\medskip

Finally an axiom for \emph{multiplying} transfinite sums:
 \item[\PA]\ (\emph{\textbf{Product axiom}}):  \ \ \ \ 
$\ \ \ (\sum_{\ag}x_{\ag})(\sum_{\bg}y_{\bg})=\sum_{\ag,\bg} x_{\ag}y_{\bg},\ \ \ $

\medskip
where the \emph{transfinite double sum}\ $\sum_{\ag,\bg}$  is defined by: 
\begin{equation*}
\sum_{\ag,\bg}\!\! x_{\ag\bg}= \!\!
\sum_{\cg} \big(\!\!\sum_{\ag\vee\bg=\cg}\!\!x_{\ag\bg}\big)
\ \  \textrm{where}\ \  \o+_{\dg\in I} 2^\dg\vee\o+_{\dg\in J}2^{\dg}=\!\!\o+_{ \dg\in I\cup J}\!\!2^\dg\ 
\end{equation*}
%
%
%
   
    (\ie\ the ordinal $\ag\vee\bg$ is the supremum w.r.t. the lattice ordering $\qincl$.)
\end{description}
    
    \medskip\noindent
     So   $\sum_{\cg=\ag\vee\bg}x_{\ag\bg}$ is an \emph{ordinary finite sums of integers}, and the corresponding comparison criterion holds:    
\begin{equation*}
\exists\,\thg\!<\kg\  \forall\, \dg\lincq \thg\,  \Big(\sum_{\ag,\bg\qincl\dg}x_{\ag\bg}\,\leq\,\sum_{\ag,\bg\qincl\dg
}y_{\ag\bg}\Big) 
  \ \Imp\  \ 
\sum_{\ag,\bg}x_{\ag\bg}\leq\sum_{\ag,\bg}y_{\ag\bg} 
\end{equation*}

Then the following useful property follows:

\medskip
\noindent
 \textbf{\emph{Translation invariance: }}\ \
If $\eta,\cg<\kg$ and  $x_\ag=0$ for $\ag\ge 2^\eta$, then
\begin{equation*}
\TI ~~~~~~~~~~~~\sum_{\ag}x_{\ag}=\sum_{\dg}y_{\dg},\  \ \mathrm{where}\
\ y_{\dg}=%
\begin{cases}
\,x_{\ag}\!\! & if\,\,\,\,\dg=2^\eta\cg+\ag \\ 
0\!\! & otherwise%
\end{cases}%
.~~~~~~~~~
\end{equation*}

In fact,
 \ $\sum_{\ag\qincl\bg} x_{\ag}=\sum_{\dg=2^\eta\cg+\ag\qincl\bg}y_{\dg},\ $\
for  $\ag<2^\eta$ and $\bg\lincq2^\eta\cg$.

\bigskip
\begin{rem} 
An  axiom much stronger than \RA, connected with the Subtraction principle of numerosities \diff, dealt with in Subsection  \ref{add},   could characterize the \emph{non-negative} Euclidean
 integers as transfinite sums of \emph{natural numbers}:
 
\medskip
\SRA (\emph{\textbf{Strong Representation Axiom}})

\smallskip\noindent
\emph{${~~}$For all $ \xk\ge 0$ in $\Zb_\kg$\ 
there\ exists\ 
$
\mathbf{n}=\la n_{\ag}\mid \ag<\kg\ra\in\,\N^{\kg}$ \st\  $\xk= \sum_{\ag} n_\ag.$
}

\medskip\noindent
The axiom {\SRA} directly yieds  the representation axiom \RA\ above, since every Euclidean integer is the \emph{difference }of two disjoint transfinite sums of {natural numbers} $\zk=\sum_{\ag\in A^+}z_\ag-\sum_{\ag\in A^-}|z_\ag|$,
\  where $A^\pm=\{\ag\mid {z_\ag}^>_<0\}$.

\medskip
More important, the full axiom \SRA\ has the interesting consequence that all Euclidean integers are actually of the simpler form $\zk=\pm\sum_{\ag}n_\ag$ for  some $\nb\in\Nak$: this fact is of great importance for the theory of numerosities.

\end{rem}

\subsection{The counting functions}\label{count}
${~~~}$
To each\
 $\kg$-sequence of integers $\ \xb\in\Zak\ $ we associate its \emph{counting function}
$\fb_\xb:\kg\to\Z$ defined by $\ \fb_\xb(\ag)=
\sum_{\bg\qincl\ag}x_\bg$.

\begin{lem}\label{period}
Every  function $\psi:\kg\to\Z$ is 
the counting function $\fb_\xb$ of some $\kg$-sequence $\xb\in\Z^\kg$.

\end{lem}
\proof
%
Given $\psi:\kg\to\Z$
 define $x_\ag$ inductively
  by putting
 $$x_0=\psi(0),\ \
   \ x_{\ag}= \psi(\ag)-\sum_{\bg\qinclp\ag}x_{\bg}.$$
Thus $\ \ \ \psi(\ag)=\sum_{\bg\qincl\ag} x_\bg$.
\qed

\bigskip
  Then we obtain the following characterization
 \begin{thm}\label{limult}
  There exist a fine ultrafilter
 $\,\U$ over $\kg$,  
  corresponding to a prime ideal $\pk$ of the ring $\Zak$,
   and an 
 isomorphism $\sg$ of the ring $\Zb_\kg$  of the Euclidean integers onto the ultrapower $\Z\ult{\kg}{\U}\cong\Zak/\pk$ that makes the  diagram \diser\
 commute:

 \medskip
\begin{center}
\begin{picture}(90,70)
   \put(0,0){\makebox(0,0){$\Zb_\kg$}}
    \put(75,0){\makebox(0,0){$\Zak\!/\pk\cong\Z\ult{\kg}{\U}$}}
   \put(0,67){\makebox(0,0){$\Zak$}}
   \put(50,39){\makebox(0,0){\diser}}
      \put(24,25){\makebox(0,0){$\pi_\pk$}}
   \put(90,67){\makebox(0,0){$\Zak$}}
   \put(28,5){\makebox(0,0){$\sg$}}
   \put(45,74){\makebox(0,0){$\fb$}}
   \put(-8,35){\makebox(0,0){$\Sg$}}
   \put(100,35){\makebox(0,0){$\pi_{\U}$}}
   \put(8,0){\vector(1,0){40}}
   \put(7,67){\vector(1,0){72}}
   \put(0,60){\vector(0,-1){50}}
   \put(90,60){\vector(0,-1){50}}
    \put(0,60){\vector(1,-1){50}}
\end{picture}
\end{center}

\medskip
 \noindent	   
    $(\fb$ maps  $ \xb\in\Zak$ to its counting function
    $\fb_{\xb}\in\Zak$,\ $\Sg$  maps $\xb\in\Zak$ to its trans\-finite sum $\Sg(\xb)\in\Zb_\kg$, while $\pi_\pk$ and $\pi_{\U}$ are the canonical projections of  $\,\Zak$ onto  the quotient ring $\!\!\!\mod\pk$ and onto the ultrapwer $\!\!\!\!\mod\U,
    respectively)$

Moreover $\Sg$ maps the semiring $\Nak$ of all  $\kg$-sequences of natural numbers onto the nonnegative part $\Zb_\kg^{\ge 0}$ of $\ \Zb_\kg$ if and only if $\  \Zb_\kg$ satisfies the strong representation axiom \SRA.

 \end{thm}
 
 \proof 
Let  $\U\linc\C_\kg$ be the ultrafilter generated by the zero-sets of $\Sg$, \ie\ the sets  $\Zl(\xb)=\{\ag<\kg\mid \fb_\xb(\ag)=0\}$ of those $\xb\in\Zak$ that give $\Sg(\xb)=0$.  Then $\sg$ is consistently and uniquely defined by putting 
$\ \ \sg(\Sg(\xb))=\pi_\U(\fb_\xb)$. 

 It turns out that any \emph{fine ultrafilter} $\U$ on $\kg$ validates the axioms \RA,\LA, \CA, and \PA. (Remark that $\U$ must intersect all cones $C(\thg)$, by axiom \CA, hence it should be fine.)
 \qed

\bigskip
  Thus the existence of rings of Euclidean integers is granted, and
    problems might arise only in the search  for suitable ultrafilters validating the strong representation axiom \SRA, which, differently from the other axioms, is independent  even of the strong axioms of the Euclidean field $\EE$ of \cite{BFeu}.
   
 \medskip
\begin{rem} 
By virtue of the fine ultrafilter $\U$, the comparison axiom may receive the following formulation
\item[~\CA] ({\emph{\textbf{Comparison Axiom}})} ${~}$
For all $\mathbf{x,y}\in{\,\Zak}$ 
\begin{equation*}
 \sum_{\ag}x_{\ag}\leq\sum_{\ag}y_{\ag}~
   \ \Iff\ ~ \ \exists\,U\in\U\  \forall \, \dg\in U\, \Big( \ \sum_{\ag\qincl \dg}x_{\ag}\,\leq\,\sum_{\ag\qincl\dg
}y_{\ag}\,\Big)
\end{equation*}

\medskip
\noindent
In particular, if the ultrafilter $\U$ is\emph{ superfine}, \ie\ contains, besides the cones $C(\theta)$, also the sets 
$\ D(\eta)=\{\dg=2^{\eta\cdot 2}\ag+2^\eta\xi+\xi\mid \xi<2^\eta,\,\ag<\kg\},$
for all $\eta<\kg$,
then   criterion \textbf{(D)} of Subsection \ref{record} gives the following property:

\medskip
\noindent
 \textbf{\emph{Double sum linearization: }} 
If $\eta<\kg$ and $x_{\bg,\cg}=0\  for \  \cg,\bg\ge2^\eta
$, then
\begin{equation*}
 \DSL { ~~~~~~~\sum_{\bg,\cg}x_{\bg\cg}=\sum_{\ag}y_{\ag},\  \ \mathrm{where}\
\ y_{\ag}=%
\begin{cases}
\,x_{\bg\cg}\!\! & if\,\,\,\,\ag=2^\eta\bg+\cg \\ 
0\!\! & otherwise%
\end{cases}%
.~~~~~~~}
\end{equation*}

\noindent
In fact,  by taking  $\dg\in D(\eta)$, 
$~\sum_{\bg,\cg\qincl\dg}x_{\bg\cg}=
\sum_{2^\eta\bg+\cg\qincl\dg}x_{\bg\cg}=\sum_{\ag\qincl\dg} y_\ag.
$\

\end{rem}
  
   \subsection{Embedding the ordinals}\label{embord}
The  fact that the ordinals less than $\kg$, with the so called ``natural sum and product" can be naturally embedded into $\Zb_\kg$ as an ordered subsemiring, could be proved exactly as in   \cite{BFeu},
but it should be remarked that, in dealing with products,  \cite{BFeu} uses the property $\DSL,$  that follows there from the axiom \DS, that is not assumed here.
However, as proved above, the property $\DSL$ holds when the ultrafilter $\U$ of Theorem \ref{limult} is superfine.
\begin{thm}\label{ordeu}{\emph{ (Thm. 2.2 of \cite{BFeu}) }}
 Define $\Psi:\kg\longrightarrow\Zb_\kg\ $  by\ $\
\Psi(\bg)=\sum_{{\alpha}}\chi_\bg(\ag) $,\\ 
\smallskip
where $\chi_\bg$ is the characteristic function of $\bg$,  \ie\ \ $\chi_\bg(\ag)=\begin{cases}
    1  & \text{ if } \ag<\bg, \\
     0 & \text{otherwise}
\end{cases}$.
 Then
$~~~~~~~~ \ag<\bg \iff \Psi(\ag)<\Psi(\bg), \ \ and\ \ \ \Psi(\alpha\oplus\beta)=\Psi(\alpha)+\Psi(\beta),\  $\\
hence $\Psi$ is an isomorphic  embedding of $(\kg;<,\oplus)$,  as ordered  semigroup, into the nonnegative part of $\Zb_\kg$.

Moreover, if the property $ \DSL$\ holds, then also
$\ \ \Psi(\alpha\otimes\beta)=\Psi(\alpha)\cdot\Psi(\beta),\ \ $
and $\Psi$ is an isomorphic  embedding of $(\kg;<,\oplus,\otimes)$,  as ordered  semiring, into the nonnegative part of the ring $\Zb_\kg$.


\rm{(Recall that  $\oplus, \otimes$  denote the \emph{natural} sum and product of ordinals.)}
\end{thm}
\proof
The values of $\Psi$ being transfinite sums of ones ``without holes", the assertion 
on $<$ is immediate. Moreover,
commutativity, associativity and distributivity holding  both in the ring $\Zb_\kg$ and in 
the semiring $(\kg,\oplus,\otimes)$, it suffices to  refer to the base-$2$ normal form and prove that,  for $\theta\ge\eta$, 
$$\Psi(2^\theta)+\Psi(2^\eta)=\Psi(2^\thg\oplus2^\eta)\ \ \rm{and}\
\  \Psi(2^\theta\otimes 2^\eta)=
\Psi(2^{\theta\oplus\eta})=\Psi(2^\theta)\cdot\Psi(2^\eta).$$
The first equality follows directly by Translation Invariance $\TI$, whereas the latter follows when the property $\DSL$
of Double Sum Linearization holds:
 $$\Psi(2^\theta)\cdot\Psi(2^\eta)=(\sum_\bg \chi_{2^\theta}(\bg))(\sum_{\cg}\chi_{2^\eta}({\cg}))=
\sum_{\bg,\cg} w_{\cg\bg},\
\ \  w_{\cg\bg}=
\begin{cases}
1\!\! & if\ \bg<2^\theta, \cg<2^\eta \\ 
0\!\! & otherwise%
\end{cases}
$$
$$\Psi(2^\theta\otimes 2^\eta)=  \Psi(2^{\theta\oplus\eta})=\sum_\ag \chi_{2^{\theta\oplus\eta}}(\ag)=
\sum_\ag z_\ag,\  \ \ z_{\ag}=
\begin{cases}
w_{\cg\bg}\!\! & if\, \ag=2^\theta\cg+\bg\\ 
0\!\! & otherwise%
\end{cases}
$$
Actually,   by taking any $\dg\in D(\thg)$, criterion \textbf{(D)}\ gives
$$~~~~~~~~~\sum_{\ag\qincl\dg} z_\ag=
\sum_{2^\eta\cg+\bg\qincl\dg}w_{\cg\dg}=\sum_{\cg,\bg\qincl\dg}w_{\cg\dg},
~$$
%
and so
$~\sum_{\ag} z_\ag=
\sum_{\cg,\bg}w_{\cg\bg}
$ \ follows by $\DSL$.
\qed

\bigskip

\medskip
It is worth noticing that the meaning of the \emph{natural product} between ordinal numbers, defined
through order types, might seem quite involved and \emph{not easily intuitive}. On
the contrary, thinking of an ordinal number as a particular \emph{Euclidean integer},
namely as a \emph{transfinite sum of ones ``without holes"}, makes appear quite natural the
meaning of the product, as given by the product formula. And by the same reason the
natural ordering of ordinals obviously agrees with that induced by the ordering of \,${\Zb_\kg}$.

 \subsection{The Euclidean ultrafilters.}\label{cons}
 
  Given a function $\psi\in\Zak$ that is positive modulo $\U$,  in order to validate 
  the strong representation axiom \SRA\ 
  one needs a nonnegative sequence $\xb\in\Nak$ such that  $\ \{\ag<\kg\mid \fb_\xb(\ag)=\sum_{\bg\qincl\ag}x_\bg=\psi(\ag)\}\in\U$.\
   Now, nonnegative $\kg$-sequences $\xb$ give rise to nondecreasing counting functions $\fb_\xb$, and conversely, so one needs  a fine ultrafilter $\U$ including, for each $\psi\in\Nak$, a set $U_\psi$ such that
 $$(\#)~~~~~~~~~~~~~~~~~\forall \ag,\bg\in U_\psi\ \big( \ag\qinclp\bg\ \ \Imp\ \ \ \psi(\ag)\le\psi(\bg)\big).~~~~~~~~~~~~~~~~$$
 
 Let $~[\kg]^2_\qinclp=\{(\ag,\bg)\mid\ag\qinclp\bg\}$ be the set of all 
   $\,\qinclp$-ordered pairs, and let  \ $G: [\kg]^2_\qinclp\to \{0,1\}$
   be a  $2$-partition of $[\kg]^2_\qinclp$:
   a set $H_i\incl\kg$ \ is\ $i$-\emph{homogeneous for} $G$ if
   $\ \forall\ag,\bg\in H_i\ (\ag\qinclp\bg\ \Imp\ G(\ag,\bg)=i).$

  For $\psi\in\Nak$, define the partition
 
 \vspace{-0,5cm}
 $$G_\psi:[\kg]^2_\qinclp \to \{0,1\}\ \ \textrm{by}\ \ G_\psi(\ag,\bg)=\begin{cases}
  0    & \text{if } \psi(\ag)>\psi(\bg), \\
   1   & \text{otherwise}.
\end{cases}$$

Then the partition $G_\psi$ does not admit 
any 
sequence $\la\ag_n\mid n<\og\ra$ \st, for $n<\og$, 
$\ \ag_n\qinclp\ag_{n+1}$\ {and}
\ $G(\ag_n,\ag_{n+1})=0$
(call such a sequence a \emph{$0$-chain}).

Call \emph{Euclidean} a fine ultrafilter $\U$ on $\kg$  if for all 
$\psi\in\Nak$ the $2$-partition
 $G_\psi$  of $[\kg]^2_\qinclp$ has a  homogeneous set $U_\psi\in\U$.
Clearly $U_\psi$ is $\qinclp$-cofinal, hence it cannot be $0$-homogeneous, so it satisfies the condition $(\#)$.
Hence Theorem \ref{limult} yields the following characterization:
\begin{cor}\label{exis}
The ring  of the Euclidean integers
 $~\Zb_\kg\cong\Z\ult{\kg}{\U}$  satisfies the strong representation axiom \SRA\ if and only if  the ultrafilter $\,\U$ is Euclidean.
 \qed
 \end{cor}

 We are left with the question of the existence of Euclidean ultrafilters.

\begin{rem}\emph{The partition property $\ [A]^{<\og}_\pincl\to(cofin)^r_k$.}\ (see \cite{Jecof,JS})\label{part}

\smallskip
The partition property $[A]^{<\og}_\pincl\to(cofin)^r_k$ means that any finite partition of all $\pincl$-ordered $r$-tuples of finite subsets of $A$ admits a $\pincl$-cofinal homogeneous subset $H\incl A$, \ie\  every $u\in [A]^{<\og}$ is included in some $v\in [H]^{<\og}$, and all $\pincl$-ordered $r$-tuples from $[H]^{<\og}$ belong to the same piece of the partition.

Clearly, the partition property $[A]^{<\og}_\pincl\to(cofin)_k^2\,$ depends only on $|A|$.
For countable $A$ it
follows immediately from Ramsey's Theorem, and for $|A|=\ah_1$ it
has been proved by Jech and Shelah in \cite{JS}, while the general problem, already posed by Thomas Jech in 1973, remains still open (see \cite{Jecof}).
\end{rem}

%
Considering the corresponding notion for the formal inclusion $\qinclp$, the validity of $\ \kg_\qinclp\to(\qinclp$-$cofin)^2_2$  would yield directly  the existence of Euclidean ultrafilters on $\kg$; 
however,  the full property $\ \kg_\qinclp\to(\qinclp$-$cofin)^2_2$  might be stronger than needed for the existence of Euclidean ultrafilters: to be sure, a  (possibly weaker) property providing
 an appropriate version of the Erd\"os-Dushnik-Miller partition property $\kg\to(\og,\kg)^2$ suffices.
 
 \subsubsection{The partition property $\ \kg\to(\og,\kg)^2_\qinclp.$}\label{part}
 \begin{defn}
The partition property  $\ \kg\to(\og,\kg)^2_\qinclp\ $  affirms that \emph{any $2$-partition  $\ G: [\kg]^2_\qinclp\to \{0,1\}\ $ either admits a $0$-chain} (\ie\  a $\qinclp$-increasing sequence $\ag_n$ with $G(\ag_n,\ag_{n+1})=0)$, or it has a $\qinclp$-cofinal homogeneous set 
$H$ (hence necessarily $1$-homogeneous). 
\end{defn}
\bigskip
This partition property is all that is needed in order to have Euclidean ultrafilters, namely
\begin{lem} \label{ultEu}
 If $\ \kg\to(\og,\kg)^2_\qinclp $ holds, then there are Euclidean ultrafilters on $\kg$.  
\end{lem}
\proof \ Given $\psi_1,\ldots,\psi_n\in\Nak$, define $G_{\psi_i}$ as above:
then the product partition
 $G_\psi=\prod_1^n G_{\psi_i}$ cannot admit $0$-chains,
  so there is a $\qinclp$-cofinal   $1$-homogeneous set $H_\psi$, which is  
   simultaneously $1$-homogeneous for all $G_{\psi_i},\ 1\le i\le n$.
 Hence the family  $\,\Hl=\{H_\psi\mid \psi\in\Nak\}$ 
  has the \FIP,
and any fine  ultrafilter $\U$ on $\kg$ including $\Hl$ is Euclidean.
\qed

\bigskip
It turns out that the above property $\ \kg\to(\og,\kg)^2_\qinclp $, has been recently stated for all cardinals $\kg$ in \cite{EDM}, where the authors
follow the track traced in section 2 of \cite{JS}, with the ordinal $\ag$ identified with the finite set $L_\ag$, and formal inclusion $\qinclp$ replacing ordinary set-inclusion $\pincl$.
So
the existence of rings of Euclidean integers satisfying the strong representation axiom \SRA\ is granted.

 \section{Euclidean measures of size for sets}\label{eusiz}
 In set theory the usual measure of the size of sets is  
 is given by the classical Cantorian
 notion of \emph{``cardinality''},  whose ground is the so 
 called \emph{Hume's Principle}
 \begin{center}
      \emph{Two sets 
  have the same size if and only if there exists a biunique 
  correspondence between  them.}
 \end{center}
This 
assumption might seem natural, and even \emph{implicit in the notion of 
counting}; 
but it strongly violates the equally natural \emph{Euclid's principle} applied to sets
 \begin{center}  
    \emph{A set 
     is greater than its proper subsets,}
 \end{center}
which in turn seems implicit in the notion of  \emph{magnitudo}, even for sets.

So one could distinguish  two basic kinds of size theories for sets:

\begin{itemize}
\item A  size theory  is  \emph{Cantorian}
if, for all $A,B$:

    $(\mathsf{{HP}})$ \ \
    ${~~~~~~~~~~~~~}\ A\eq B\ \ \Iff\ \ \exists f:A\to B$ biunique
\item A  size theory is  \emph{Euclidean}
if, for all $A,B$: 

   $(\mathsf{{EP}})$ \ \
    ${~~~~~~~~~~~~}\ \   \ A\< B \ \Iff\ \ \exists C\ s.t.\ A\pincl C\ \&\ B\eq C.$\\
(Remark the use of \emph{proper inclusion} in defining \emph{strict comparison} of sets.)
\end{itemize}
The consistency of the principle \EP\ for  \emph{uncountable sets} appeared problematic from the beginning, and this question has been posed in several papers (see \cite{BDNFar,BDNFuniv,DNFtup}), where only
the \emph{literal set-theoretic translation} of
      the fifth Euclidean notion, \ie\  the sole left pointing arrow of $ \EP$,  
    $$ \Ec{~~~~~~~~~~~~~~~~~~~~~~~~~~  ~}A\pincl B\ \ \Imp\ \ A\< B,{~~~~~~~~~~~~~~~~~~~~~~~~~~~~}$$
%
%
       has been obtained,
 (On the other hand, it is worth recalling that also the totality of the Cantorian weak cardinal ordering 
had to wait more than two decades till \emph{Zermelo's new axiom of choice} to be established!)

A general discussion of different ways for comparing and measuring the size of sets can be found in \cite{Compar}. Here we present 
 a Euclidean numerosity theory for suitable collections $\WW$ of point sets of  finite dimensional spaces over lines  $\LL$ of arbitrary cardinality,  satisfying the \emph{full principle} \EP.  This numerosity might be extended to the whole universe under 
simple  set theoretic assumptions, \pes\ Von Neumann's axiom, that gives a (class-)bijection between the universe $V$ and the class $Ord$ of all ordinals.
 
 \bigskip
 \subsection{Natural congruences}\label{cong}
First of all, once the general Hume's principle cannot be assumed,
 the fourth 
Euclid's common notion
\begin{center}
    \emph{Things exactly applying  onto one 
     another are equal to one another}
\end{center}
is left in need of an adequate choice of natural ``{exact applications" that preserve size} ($=$\emph{congruences}):
so we  isolate the group $\Gk(W)$ of  the ``natural transformations'' 
 of tuples, \ie\  those 
 preserving the \emph{support} (the set of components) of a tuple
   $(supp(a_1,\ldots,a_n)=\{a_1,\ldots,a_n\})$\footnote{
~in particular
$A\times B\eq B\times A$ and 
 $(A\times B)\* C \eq A\* (B\times C)$.},   that seem appropriate for a Euclidean theory involving sets of tuples,
and postulate 
\begin{itemize}
  \item[\CP] (Congruence Principle)
   $\ \sg\in\Gk(W)\  \Imp\  \forall A\in \mathbb{W}\  \big(\sg[A]\in\WW \ \mathrm{\&} \ \sg[A]\simeq A\big).$
\end{itemize}

%

\subsection{Addition of numerosities}\label{add}
One wants not only \emph{compare}, but also \emph{add and subtract
 magnitudines}, according to the second and third Euclidean common notions

\emph{{~~~~~~~~}... if equals be added to equals, the wholes are equal.{~~~~~~~~~~~~~~~~}}

\emph{...  if equals be subtracted from equals, the remainders are equal.}

\medskip
When dealing with sets, it is natural to take \emph{addition} to be \emph{(disjoint) union}, and \emph{subtraction} to be \emph{(relative) complement}, and
\smallskip  \noindent  
 the following 
\begin{description}
\item[$(\mathsf{{AP}})$]  (\emph{Aristotle's Principle})\footnote{~
This priciple has been named \emph{Aristotle's Principle} in \cite{QSU,FM}, 
because it resembles     Aristotle's preferred example of a ``general axiom". It is especially relevant in this context, because \AP\ implies both the second and the third 
     Euclidean common notions, and also the fifth 
     whenever no 
     nonempty set is equivalent to $\0$, as stated in the proposition below. 
}
    ${~~~~~~~~~~~}
    A\eq B\ \Longleftrightarrow\
    A\7 B \eq B\7 A.
    $
\end{description}
is convenient, because it yields both the second and third Euclidean common notions, 
see \cite{QSU}.

 Define the
 addition  of \emph{numerosities}, \ie\  equivalence classes modulo $\eq$,
by 
 $$[A]+[B]=[A\cup B]\ \ for\ all \ \ A,B \  such\ that\ \ A\cap B=\0;$$
then the quotient set $\Nk=\WW/\eq$  becomes a \emph{positive semigroup},
\ie the non-negative part of an ordered abelian group.
In particular one gets the  \emph{``most wanted Subtraction Principle"} of \cite{BDNFar}:\\
\ \diff~~~~~~~~~~    \emph{ $A\< B\ \ \Iff\ \ \exists C\ne\0  \ \big(C\cap (A\cup B)=\0,
      \  
    (C\cup A)\eq B\Big)$.~~~~~~~~~~~~~~~~~~}\label{dif}

The  consistency problem of the Subtraction Principle, 
studied in several papers dealing with Euclidean (also called \emph{Aristotelian}) notions of size for sets, receeved a positive answer only for \emph{countable} sets in \cite{DNFtup,QSU,FM}. A positive answer for sets of arbitrary cardinality is obtained in \cite{Compar}, and follows from the existence of Euclidean ultrafilters, see Subsubsection \ref{numeu}).

\subsection{Multiplication of numerosities}\label{mult}
In classical mathematics,  geometric figures having different dimensions are never compared,  so a \emph{multiplicative} version of Euclid's second common notion
\begin{center}
   \emph{...  if 
        equals be multiplied by equals, the products are equal}
\end{center}
was not considered, for ``{dishomogeneous}" magnitudes.

On the other hand,  in modern mathematics a single
set of ``numbers", the real numbers $\R$,
is used as a common scale for
  the size of figures of any dimension.
In a general set theoretic context it seems natural to consider abstract sets as homogeneous mathematical objects, without distinctions based on dimension,
%
and a satisfying \emph{arithmetic} of numerosities needs a \emph{product} (with a corresponding \emph{unit}):
we adhere to the natural Cantorian choice
 of introducing multiplication 
through \emph{Cartesian products},\footnote{~
 Although the Cartesian product is neither commutative nor associative \emph{stricto sensu}, nevertheless the corresponding natural transformations should be taken  among the congruences in the group $\Gk(W)$.
}
 and  taking \emph{singletons as unitary.}\footnote{
\emph{CAVEAT}: the Cartesian product is optimal when any two sets  $A, B$ are \emph{multipliable} in the sense that their Cartesian product is disjoint from their union, but
 when \emph{transitive universes} like $V_\kg, H(\kg)$, or $L$ are considered becomes untenable (\pes\ already $V_\og\*\{x\}\pincl V_\og$ for any $x\in V_\og$),(see \pes\  the discussion in \cite{BDNFar,Compar}), hence
 not all singletons may be ``suitable" for a Euclidean theory.}
So any $A\*\{b\}, b\in B$ may be viewed as a \emph{disjoint equinumerous copy of }$A$, thus  making 
  their (disjoint) union $A\*B$  the sum of ``$B$-many copies of $A$", in accord with the intuitive idea of product.    

\subsection{Euclidean numerosities as Euclidean integers} \label{numeu} 
 The properties of the ring $\Zb_\kg$ of the Euclidean integers allow 
 for a simple engrafting of a Euclidean numerosity  for ``Punktmengen",
  \ie\ sets of tuples, over any line $\LL$ of arbitrary size $\kg$.
Since in a general set-theoretic context there are no ``geometric" or ``analytic" properties to be considered, the sole relevant characteristic of the line $\LL$ remains  cardinality, so a convenient choice seems to be simply identifying $\LL$ with its cardinal $\kg$, thus obtaining the fringe benefit that no pair of ordinals is an ordinal, and Cartesian products may be freely used.
Grounding on the preceding discussion, we  pose the following definition
\begin{defn}\label{eunum}
A 
 \emph{(Euclidean) numerosity} for ``Punktmengen" over $\LL$ is a pair $(\WW,\ll)$, where $\ll$ is a total preordering on a set $\WW\incl \P(\bigcup_{n\in\N} \LL^n)$ \st\
\begin{center}
  $A\cup B,\ A\*B\in\WW\ \ \Iff\  A,B\in\WW,\ \ \  C\incl A\in\WW\ \Imp\ C\in\WW$
  \end{center}
and  the following conditions are satisfied for all $A,B,C\in\WW$:
\begin{itemize}
 \item[$(\mathsf{{EP}})$] $ A\prec B\ \Iff\ \exists B'(A\subset B'\eq B),$
where  \emph{inclusion} and \emph{preordering} are \emph{strict};
\item[\CP]  $\tau[A]\eq A$ for all $\tau \in\Gk(\LL),$ the group of all \emph{support-preserving} bijections;\footnote{~see Subsubsection \ref{cong}; in particular
 $A\times B\eq B\times A$ and 
 $(A\times B)\* C \eq A\* (B\times C)$, so commutativity and associativity of multiplication follow.
}

 \item[$(\mathsf{{AP}})$]$A\eq B\ \Longleftrightarrow\ A\7 B \eq B\7 A;$
  \item[\PP] 
 $A\eq B \ \ \Iff \ \ A\* C\eq B\* C$ (for all $C\ne\0$); 
\item[\UP] $A\eq A\* \{w\}$ for all $w\in W=\bigcup\WW$.
\end{itemize}
\noindent
\end{defn}
 The above  Principles provide the \emph{set of numerosities}  $\Nk=\mathbb{W}/\eq$
with the best arithmetic properties, namely those of the nonnegative part of an ordered domain, see \pes\ \cite{Compar}. 
%
%
 Here, having at disposal  the ring $\Zb_\kg$ of the Euclidean integers,
we give first a  function 
$\nk:\P(\kg)\to \Zb_\kg$ as the transfinite sum of the characteristic functions of 
each subset of $\kg$:
$$\nk(A)=\sum_{\ag} \chi_{A}(\ag),\ \ where\ \ \chi_A(\ag)=\begin{cases}
    1  & \text{ if\ } \ag\in A\\
   0   & \text{otherwise}.
\end{cases}
  $$
 Thus in particular $\nk(\ag)=\Psi(\ag)$ for all ordinals. 
 
Then we may extend the function
$\nk$ to $\P(\kg^n) $ by
 assigning to each $n$-tuple 
$\_\ag=(\ag_{1},\ldots,\ag_{n})\in\kg^n$ the ordinal 
$ \psi_n(\_\ag)=\ag_1\lor\ldots\lor\ag_n<\kg$, and putting,
for $A\incl\kg^n$, 
$\ \nk(A)=\sum_{\ag} \chi^{(n)}_{A}(\ag),\ \ \textrm{where}\ \
 \chi^{(n)}_A(\ag)=|\{\_\ag\!\in\! A\mid \psi_n(\_\ag)=\ag\}| $

Remark that we are assigning the same ordinal $\ag$ to $\ag\in\kg$, to $(\ag_1,\ag_2)\in\kg^2$  if $\ag=\ag_1\lor\ag_2$, \ldots, to $(\ag_{1},\ldots,\ag_{n})\in\kg^n$   if $\ag=\bigvee_1^n \ag_i$,  hence the functions
 $\chi^{(n)}$, for $n>1$, are not properly characteristic functions, but they assume nonnegative integer values, so their sums are nonnegative Euclidean integers. 

Now we can easily extend the numerosity function $\nk$ to all \emph{finite dimensional} point sets, \ie\ sets $A$ \st\ $\{n\mid A\cap\kg^n\ne\0\}$ is finite, namely
$$\nk(A)=\sum_n \nk(A\cap\kg^n).$$ 
Clearly \begin{itemize}
  \item $\ A\pincl B\ \Imp\ \nk(A)<\nk(B), $\  so the fifth Euclidean notion  \Ec\ holds;
  \item the Aristotelian principle \AP\ holds because
$\nk(A)=\nk(A\7 B)+\nk(A\cap B)$;
  \item the principles  \UP\ and \PP\ follow directly from the product axiom \PA.
\end{itemize}
Finally, the principle \EP\ is here equivalent to  the strong representation axiom \SRA, because both postulate that the difference $\chi_A-\chi_B$  of two  (generalized) characteristic functions, when positive,  has the same transfinite sum of a single positive function $\chi_C$.

Summing up, $\Nk=\Zb_\kg^{\ge 0}$ is a positive semiring where the five Euclid's Common Notions are satisfied,  all finite sets receive their number of
elements as numerosity,  hence it contains as initial segment an isomorphic copy of 
the natural numbers $\N$; moreover one obtains the supplementary benefits that every point set is equinumerous to a set of ordinals,  and conversely  that every nonnegative Euclidean integer $\xk$ is the numerosity of a set $X$ of ordinals, namely the transfinite sum of the characteristic function $\chi_X$.

\section{Final remarks and open questions}\label{froq}

\subsection{The Weak Hume Principle and the Subset Property}\label{hum}

Perhaps the best way to view a Euclidean numerosity is looking at it as a \emph{refinement} of Cantorian cardinality, able to separate sets that, although \emph{equipotent}, should have in fact \emph{really different sizes}, in particular when they are \emph{proper} subsets or supersets of one another. To this aim, 
  the principle \EP\  might be integrated by adding the clause
 $$(\mathsf{{WHP}}){~~~~~~~~~~~~~~~~~~~~~~~~} \nk(A)\le\nk(B)\ \ \Imp\ \ \ \exists 
  f\, 1\mbox{- to -}1,\, f: A\to B.{~~~~~~~~~~~~~~~~~~~~~~~~~~}$$

 So sums of ones of greater cardinality produce greater Euclidean integers, hence the ordering of the Euclidean numerosities refines the cardinal ordering, satisfying  the ``Weak Hume's Principle" 
  \begin{center}
      \ \emph{If two sets 
 are equinumerous, then there exists \\ a biunique 
  correspon\-dence between  them.}
 \end{center}
 
Another interesting consequence of $(\mathsf{{WHP}})$ is the property 

\smallskip\noindent
$\oob$  (\emph{Proper Subset Property})
 $~\ \ ~A\prec B\ \ \Iff\ \ \exists A'\pincl B\ \ \textit{s.t.}\ \  A'\eq A.$
 
\medskip\noindent
In general, the set of numerosities  has the same size as the universe $W$,
  since one can define strictly increasing chains of sets of arbitrary length; but any  set $A$ has only $2^{|A|}$ subsets, and so the
Proper Subset Property  \oob\ implies that
 \emph{the initial segment of numerosities}  generated by $\nk(A)$ has size $2^{|A|}$, 
contrary, \pes, to the  large ultrapower  models of \cite{BDNFar,BDNFuniv}.

This topic is dealt with in \cite{Compar}, where it is proved that
the family of sets  
 $$Q^{>}_{AB}=\{\bg<\kg\mid \sum_{\ag\qincl\bg}\chi_A(\ag)> \sum_{\ag\qincl\bg}\chi_B(\ag)\}\ \textrm{for }\  |A|>|B|.$$
has the \FIP\ together with the cones $C(\thg)$, hence may be contained in the fine ultrafilter $\U$. However it is not known whether that family might be included in an Euclidean ultrafilter, so the consistency of the weak Hume principle \WHP\ with the difference property \diff\ is still open.

\subsection{The power of numerosities}\label{pow}
The power $\mk^{\nk}$ of infinite numerosities  is here always well-defined, since 
numerosities are
    \emph{positive euclidean numbers}, hence \emph{nonstandard natural numbers.}
By using finite  approximations given by  intersections with \emph{suitable finite sets,}     the interesting  relation
    $$2^{\nk(X)} = \nk([X]^{<\og}),$$
    has been obtained already  in \cite{BDNFar}.
Since the comparison axiom \CA\ evaluate transfinite sums by finite sums $\sum_{\ag\qincl\dg}\chi_A(\ag)$, one obtains  the following general set theoretic interpretation of powers:
   $$\mk(Y)^{\nk(X)} = \nk(\{f:X\to Y\mid |f| <\aho\}),$$ 
 by considering sets of ordinals $X,Y$, and labelling each finite function $f$ by the $\qincl$-supremum of the (finitely many) ordinals involved in $f$.
    
The interesting problem of finding  appropriately defined arithmetic operations that give instead the numerosity of
    the \emph{full powersets} and \emph{function spaces} requires a quite different
    approach, and the history of the same problem for cardinalities suggests that it could not be properly solved.

\end{document}